\documentclass{article}

\usepackage{amssymb}
\usepackage{amsmath}
\usepackage{times}
\usepackage{graphicx}

\numberwithin{equation}{section}

\def\FermatC{{C}}
\def\u{{u}}
\def\v{{v}}
\def\y{{y}}

\def\squareforqed{\hbox{\rlap{$\sqcap$}$\sqcup$}}
\def\qed{\ifmmode\squareforqed\else{\unskip\nobreak\hfil
\penalty50\hskip1em\null\nobreak\hfil\squareforqed
\parfillskip=0pt\finalhyphendemerits=0\endgraf}\fi}

\newtheorem{theorem}{Theorem}{\bfseries}{\itshape}
\newtheorem{lemma}{Lemma}{\bfseries}{\itshape}
{\bfseries}{\itshape}
\newtheorem{proposition}{Proposition}{\bfseries}{\itshape}
\newtheorem{remark}{Remark}{\bfseries}{\itshape}
\newtheorem{definition}{Definition}{\bfseries}{\itshape}
\newtheorem{example}{Example}{\bfseries}{\itshape}

\newenvironment{proof}{\par\noindent{\bf Proof.}}{}

\begin{document}

\title{Decomposed Richelot isogenies of \\ Jacobian varieties of hyperelliptic curves \\ and generalized Howe curves
}

\author{Toshiyuki Katsura  \\ The University of Tokyo, \ \ \texttt{tkatsura@g.ecc.u-tokyo.ac.jp}           %
\and
        Katsuyuki Takashima \\ Waseda University, \ \ \texttt{ktakashima@waseda.jp}
}


\maketitle

\begin{abstract}
We advance previous studies on decomposed Richelot isogenies (Katsura--Takashima (ANTS 2020) and Katsura (J. Algebra)) which are useful for analysing superspecial Richelot isogeny graphs in cryptography.  
We first give a characterization of decomposed Richelot isogenies between Jacobian varieties of hyperelliptic curves of any genus. We then define generalized Howe curves, and present two theorems on their relationships with decomposed Richelot isogenies.  We also give new examples including a non-hyperelliptic (resp.\,hyperelliptic) generalized Howe curve of genus 5 (resp.\,of genus 4).
%
\end{abstract}


\section{Introduction}


\subsection{Background}
\label{Sec:Background}

Richelot isogenies of Jacobian varieties of nonsingular projective curves are generalizations of 2-isogenies of elliptic curves (see Definition \ref{Def:Richelot} for the detail),
and such isogenies of Jacobians of superspecial genus-1 and 2 curves are frequently used in post-quantum cryptography, which remains secure even when large scale quantum computers are deployed for cryptanalysis. Consequently, intensive study  (\cite{ChaGorLau09,Tak17,FT19,CDS,CS20} etc.) has been devoted to security evaluation of the isogeny-based cryptography in recent years. Here, cryptographic operations consist of random walks on graphs of isogenies between Jacobian varieties of superspecial curves. 

Costello and Smith \cite{CS20} used ``decomposed'' subgraphs of the superspecial isogeny graphs for their cryptanalysis successfully, in which {\em decomposed} principally polarized abelian varieties are cleverly used for efficiency improvements of the cryptanalysis. Richelot isogenies with such decomposed ones as codomain are called {\em decomposed Richelot isogenies}. Recent works have clarified the detailed information on such decomposed isogenies and the associated isogeny graph structures (Katsura--Takashima \cite{KT20}, Florit--Smith \cite{FS21b,FS21a}, and Jordan--Zaytman \cite{JZ}), which can be useful for more accurate analysis of the Costello--Smith attack. (See \cite{SCF22,DK23} also.) 

For a (hyperelliptic) curve $C$ of genus 2, Katsura and Takashima \cite{KT20} showed that the equivalence of existence of a decomposed Richelot isogeny from its Jacobian variety and existence of an order-2 (long) element in the reduced automorphism group. A similar equivalence for hyperelliptic curves of genus 3 was also shown by Katsura \cite{K}. These results give a basis for our present work.

Howe \cite{How16} investigated the nonsingular projective model of the fibre product of two elliptic curves (which satisfy some condition) w.r.t.\,the two hyperelliptic structures. Such curves were called Howe curves in subsequent works \cite{KHS20,KHH20}. For a genus-3 curve $C$, Katsura \cite{K} established another interesting equivalence that $C$ is a Howe curve if and only if it has a {\em completely} decomposed Richelot isogeny, whose target is a product of three elliptic curves. 
It indicates an initimate relationship between Howe curves and decomposed Richelot isogenies.  
We will further study the relationship for {\em generalized Howe curves of any genus}  (which we will define in Section \ref{Sec:Howe}).


\subsection{Our contributions}

We generalize the works \cite{KT20,K} in higher genus cases, and give a unified 
approach for investigating close connections among the three notions of decomposed Richelot isogenies, non-inversion automorphisms of order 2, and generalized Howe curves.

\begin{enumerate}

\item
We first give a decomposition criterion (Theorem \ref{Thm:hyperelliptic}) by using non-inversion automorphisms of order 2 for {\em hyperelliptic curves of any genus}, which are important in almost all cryptographic applications. As is already mentioned above, it is useful for analysing Richelot isogeny graphs in higher genus cases, where there exist works for the genus-2 case by Katsura--Takashima \cite{KT20} and Florit--Smith \cite{FS21b,FS21a}, and  for the genus-3 case by Howe--Lepr\'evost--Poonen \cite{HLP} and Katsura \cite{K}.  

\item
We then define a {\em generalized Howe curve} by the nonsingular projective model of the fibre product of two hyperelliptic curves (which satisfy some condition) w.r.t.\,the two hyperelliptic structures. We show a criterion of when a generalized Howe curve of genus $g \ge 4$ is hyperelliptic (Theorem \ref{rationalHowe}). It is simply described by using branch points of the underlying two hyperelliptic curves. As a collorary, we show that any hyperelliptic curve with an automorphism of order 2 (which is not the inverse) is realized as a generalized Howe curve (Remark \ref{Rmk:hyperelliptic}).

\item
Thirdly, we give a {\em strong} decomposition theorem for {\em generalized Howe curves of any genus} (Theorem \ref{Thm:Howe}).
Our present result is a generalization of two preceding facts: one is a complete decomposition theorem of genus-3 Howe curves \cite[Theorem 6.2]{K}, and the other is Theorem \ref{Thm:hyperelliptic} since hyperelliptic curves with a non-inversion order-2 automorphism are given by generalized Howe curves as indicated above. 

\item
We show several examples in Section \ref{Sec:Examples}. In particular, we give a generalized Howe curve of genus 5 which is non-hyperelliptic in Example \ref{Ex:NonHyperellipticLargeGenus} and that of genus 4 which is hyperelliptic in Example \ref{Ex:HyperellipticLargeGenus}, both of which are newly obtained from our theorems. 

\end{enumerate}

While we believe that our theorems, in particular, Theorem \ref{rationalHowe}, represent a meaningful advance in this research area, our understanding on the relationship between decomposed Richelot isogenies and generalized Howe curves is still slightly limited. See comments after Remark \ref{Rmk:hyperelliptic} in Section \ref{Sec:Howe}.


\subsection{Related Works}

\paragraph{Isogeny-based Cryptography}

As we already pointed out in Section \ref{Sec:Background}, a detailed study of decomposed isogenies leads to a better understanding of security of superspecial isogeny based cryptography via the Costello--Smith attack. 
In fact, very recently, Santos--Costello--Frengley \cite{SCF22} proposed to use a novel search algorithm whether a genus-2 curve has an $(N,N)$-decomposed isogenous neighbor for $2 \le N \le 11$ for improving the Costello--Smith attack. The main ingredient of their attack is given by explicit parametrizations of moduli spaces of genus-2 curves whose Jacobians have an $(N,N)$-decomposed isogeny, which are described by Kumar \cite{Kum15}. The useful, explicit descriptions depend on the special situation of genus 2. In general, it seems to be difficult to give such explicit equations of the moduli spaces for higher genera (and to efficiently compute on them). Therefore, we think that our results are a first step for employing "decomposed neighbors" for efficiently solving the higher genus isogeny problem.

As a remarkable recent progress related to decomposed Richelot isogenies, Castryck--Decru \cite{CD23} proposed a clever use of genus-2 Richelot isogenies to attack the ``elliptic curve" based SIDH key exchange protocol \cite{FJP14}. One of their key observations is that {\em decomposition events of Richelot isogenies} can be used to check right guesses among several possibilities for solving SIDH-type isogeny problems. Their attack was extended and improved by several authors soon \cite{MMPPW23,Rob23}. And, in particular, Robert \cite{Rob23} employed 8-dimensional abelian varieties and their isogenies for establishing a {\em polynomial}-time attack against SIDH protocols with {\em arbitrary} starting elliptic curves. We note that the attacks can be applied to only special cases of isogeny problems {\em with auxiliary points} as in the SIDH case, but not be applicable to the general elliptic curve isogeny problems (without auxiliary points).  

\paragraph{Enumeration of Superspecial Howe Curves}

In a subsequent work to ours, Moriya--Kudo \cite{MK22} explicitly wrote down our constructions, and established efficient algorithms for computing decomposed Richelot isogenies and generalized Howe curves in the genus-3 case. In a series of papers, Kudo and Harashita have investigated the existence and counting of superspecial curves with coauthors (see \cite{KH22} for a survey of their works). Then, Moriya--Kudo also applied their explicit algorithms to search and enumerate superspecial generalized Howe curves. We can find their Magma codes for the computations at \cite{MK22b}. 
\\

\noindent
Our paper is organized as follows: Section \ref{Sec:Preliminaries} gives mathematical preliminary results which are also shown in \cite{K,BL}. Section \ref{Sec:Hyperelliptic} presents a criterion for decomposed Richelot isogenies in the hyperelliptic curve case (Theorem \ref{Thm:hyperelliptic}). Section \ref{Sec:Howe} first defines generalized Howe curves and then gives two theorems (Theorems \ref{rationalHowe} and \ref{Thm:Howe}) on decomposed Richelot isogenies in the generalized Howe curve case. Section \ref{Sec:Examples} shows several examples.

\paragraph{Notation and conventions}
%
For an abelian variety $A$ and divisors $D$, $D'$ on $A$, we use the following (standard) notation:
\ $id_A$ and $\iota_A$ denote the identity of $A$ and the inversion of $A$, respectively.
$\hat{A}= {\rm Pic}^0(A)$ denotes the dual (Picard variety) of $A$. 
$D\approx D'$ denotes 
algebraic equivalence. 
For a vector space $V$ and a group $G$ which acts on $V$, we denote by $V^G$ 
the invariant subspace of $V$.

\paragraph{Acknowledgement} Research of the first author is partially supported by JSPS Grant-in-Aid for Scientific Research (C) No.23K03066.
Research of the second author is supported by JSPS Grant-in-Aid for Scientific Research (C) JP22K11912, JST CREST JPMJCR2113, and MEXT Quantum Leap Flagship Program (MEXT Q-LEAP) JPMXS0120319794.


\section{Preliminaries}
\label{Sec:Preliminaries}

We will review necessary mathematical preliminaries that are developed in \cite{K,BL}. 

Let $k$ be an algebraically closed field of characteristic $p > 2$.
In this section, we prepare some notation and some known lemmas 
to examine the structure of Richelot isogenies.
For an abelian variety $A$ and a divisor $D$ on $A$,
we have a homomorphism
$$
\begin{array}{rccc}
\Phi_{D} :& A &\longrightarrow & {\rm Pic}^0(A) = \hat{A}\\
     & x   &\mapsto & T_x^{*}D - D
\end{array}
$$
(cf.\,Mumford \cite{M}). Here, $T_x$ is the translation by $x \in A$.
We know that $\Phi_{D}$ is an isogeny if $D$ is ample.
    
Let $C$ be a nonsingular projective curve 
of genus $g$ defined over $k$. We denote by $J(C)$
the Jacobian variety of $C$, and by $\Theta$ the principal
polarization on $J(C)$ given by $C$.
We have a natural immersion (up to translation)
$$
\alpha_C  : C \hookrightarrow J(C)= {\rm Pic}^0(C).
$$
By the abuse of terminology, we sometimes denote $\alpha_C(C)$ by $C$.
The morphism $\alpha_C$ induces a homomorphism
$$
\alpha_C^{*} : \widehat{J}(C) = {\rm Pic}^0(J(C))\longrightarrow {\rm Pic}^0(C)=J(C).
$$
First, we give two lemmas which make clear relations of some homomorphisms.
\begin{lemma}\label{Theta}
   $\alpha_C^{*} = - \Phi_{\Theta}^{-1}$.
\end{lemma}
For the proof, see Birkenhake--Lange \cite[Proposition 11.3.5]{BL}.

Let $f : C \longrightarrow C'$ be a morphism of degree $2$
from $C$ to a nonsingular projective curve $C'$ of genus $g'$. We denote by
$J(C')$ the Jacobian variety of $C'$, and by $\Theta'$ the principal
polarization on $J(C')$ given by $C'$.
For an invertible sheaf  
${\mathcal O}_C(\sum m_iP_i) \in J(C)$ ($P_i\in C$, $m_i \in {\bf Z}$),
the homomorphism $N_f: J(C) \longrightarrow J(C')$ is defined by 
$$
N_f({\mathcal O}_C(\sum m_iP_i)) = {\mathcal O}_{C'}(\sum m_if(P_i)).
$$
By suitable choices of $\alpha_C $ and $\alpha_{C'}$, we may assume 
$$
   N_f \circ \alpha_C = \alpha_{C'}\circ f,
$$
that is, we have a commutative diagram
$$
\begin{array}{ccc}
 C &\stackrel{\alpha_C}{\hookrightarrow} & J(C) \\
            f \downarrow \quad & & \quad \downarrow N_f \\
 C' & \stackrel{\alpha_{C'}}{\hookrightarrow} & J(C').  
\end{array}
$$ 
\begin{lemma}\label{N}
  $\Phi_{\Theta} \circ f^{*} = \hat{N}_f \circ\Phi_{\Theta'}$.
\end{lemma}
For the proof, see Birkenhake--Lange \cite[equation 11.4(2)]{BL} or Katsura \cite[Lemma 2.2]{K}.

Using Lemmas \ref{Theta} and \ref{N}, we have the following lemma which is essential
to show the existence of a decomposed Richelot isogeny.

\begin{lemma}\label{2Theta}
$(f^{*})^*(\Theta) \approx 2\Theta'$.
\end{lemma}
For the proof, see Birkenhake--Lange \cite[Lemma 12.3.1]{BL} or Katsura \cite[Lemma 2.3]{K}.

\begin{definition}
Let $A_i$  be abelian varieties
with principal polarizations $\Theta_i$ ($i = 1, 2, \ldots, n$), respectively.
The product $(A_1, \Theta_1)\times (A_2, \Theta_2)\times \ldots \times (A_n, \Theta_n)$
means the principally polarized abelian variety $A_1 \times A_2 \times \ldots \times A_n$
with principal polarization 
$$
\Theta_1\times A_2 \times A_3 \times \ldots \times A_n
+ A_1 \times\Theta_2\times  A_3 \times \ldots \times A_n + \ldots +  A_1 \times A_2 \times \ldots \times A_{n-1} \times \Theta_n.
$$
\end{definition}

\begin{lemma}\label{decompose} 
Let $A$, $A_1$ and $A_2$ be abelian varieties, and let
$f : A_1 \times A_2 \longrightarrow A $ be an isogeny. 
Let $\sigma$ be an automorphism of $A$ such that 
$\sigma \circ f = f \circ(id_{A_1} \times \iota_{A_2})$
and $\Theta$ be a polarization of $A$ such that 
$\sigma^*\Theta \approx \Theta$. Then,
$$
 (A_1 \times A_2, f^*\Theta) 
 \cong (A_1, f\vert^*_{A_1}\Theta)\times (A_2, f\vert^*_{A_2}\Theta).
$$
\end{lemma}
For the proof, see Katsura \cite[Lemma 3.3]{K}.

\begin{definition}
[Richelot isogenies in genus $g$]
\label{Def:Richelot}
Let $C$ be a nonsingular projective curve of genus $g$,
and $J(C)$ be the Jacobian variety of $C$. We denote by $\Theta$
the canonical principal polarization of $J(C)$.
Let $A$ be a $g$-dimensional abelian variety with principal polarization $D$,
and $f : J(C) \longrightarrow A$ be an isogeny.
The isogeny $f$ is called a {\em Richelot isogeny} 
if $2\Theta \approx f^*(D)$. A Richelot isogeny $f$ is
said to be {\em decomposed} if there exist two principally polarized
abelian varieties $(A_1, \Theta_1)$ and $(A_2, \Theta_2)$
such that $(A, D) \cong (A_1, \Theta_1) \times (A_2, \Theta_2)$.
Moreover, the isogeny $f$ is said to be {\em completely decomposed}
if $A$ with principal polarization $D$ is decomposed into
$g$ principally polarized elliptic curves.
\end{definition}


\section{Hyperelliptic curves}
\label{Sec:Hyperelliptic}

Let $\iota$
be the hyperelliptic inversion of a hyperelliptic curve $C$ of genus $g$ ($g \geq 2$) and 
$\sigma$ be an automorphism of order 2 of $C$
which is not the inversion. We set $\tau = \sigma\circ \iota$.
We have a morphism $\psi : C \longrightarrow {\bf P}^1 \cong C/\langle \iota \rangle$.
Since the morphism $\psi$ is given by ${\rm H}^{0}(C, \Omega_C^1)$
and $\sigma$ acts on it,
the automorphism $\sigma$ induces an automorphism of ${\bf P}^1$.
In case $\sigma$ has a fixed point in the branch points of
$\psi$, $\sigma$ has two fixed points in the branch points.
Moreover, by a suitable choice of the coordinate $x$ 
of ${\bf A}^1 \subset {\bf P}^1$
we may assume that the two fixed points are given by $x = 0$ and $\infty$,
and that
$$
\sigma :x \mapsto -x; \quad y \mapsto y.
$$
Then the branch points are given by
$$
0, 1, -1, \sqrt{a_1}, -\sqrt{a_1}, \sqrt{a_2}, -\sqrt{a_2}, \ldots, 
\sqrt{a_{g-1}}, -\sqrt{a_{g-1}}, \infty.
$$
Here, $a_i$'s are mutually different and they are neither 0 nor 1. 
The normal form of the curve $C$ is given by
\begin{eqnarray*}
    y^2 = x(x^2 - 1) (x^2 - a_1)(x^2 - a_2)\ldots (x^2 - a_{g-1}).
\end{eqnarray*}
Therefore, on the curve $C$ the action of $\sigma$ is
given by
$$
    x\mapsto -x,~ y \mapsto \pm \sqrt{-1} y,
$$
which is of order 4, a contradiction. 
Therefore, $\sigma$ has no fixed points on the branch points.

Now, let the branch points be given by
$$
1, -1, \sqrt{a_1}, -\sqrt{a_1}, \sqrt{a_2}, -\sqrt{a_2}, \ldots, 
\sqrt{a_{g}}, -\sqrt{a_{g}}.
$$
Here, $a_i$'s are mutually different and they are neither 0 nor 1. 
The normal form of the curve $C$ is given by
\begin{eqnarray}
\label{Eq:NormalForm}
y^2 = (x^2 - 1) (x^2 - a_1)(x^2 - a_2)\ldots (x^2 - a_{g}).
\end{eqnarray}
Elements $x^2$ and $y$ are invariant under $\sigma$.
We set $\u = x^2$ and $\v = y$. Then,
the defining equation of the curve $C/\langle \sigma \rangle$
is given by
$$
  \v^2 = (\u - 1) (\u - a_1)(\u - a_2)\ldots (\u - a_{g}).
$$
We set $C_{\sigma} = C/\langle \sigma \rangle$.
We have the quotient morphism $f_1: C \longrightarrow C_{\sigma}$.
Elements $xy$ and $x^2$ are invariant under $\tau$. We set
$\u = x^2$ and $\v = xy$. Then,
the defining equation of the curve $C/\langle \tau \rangle$
is given by
$$
  \v^2 = \u(\u - 1) (\u - a_1)(\u - a_2)\cdots (\u - a_{g}).
$$
We set $C_{\tau} = C/\langle \tau \rangle$.
We have the quotient morphism $f_2 : C \longrightarrow C_{\tau}$.
We denote by $g(C)$ (resp.~$g(C_{\sigma})$, resp.~$g(C_{\tau})$)
the genus of $C$ (resp.~$C_{\sigma}$, resp.~$C_{\tau}$).
It is easy to see that $g=g(C) = g(C_{\sigma}) + g(C_{\tau})$.
We have a morphism
$$
   f = (f_1, f_2) :  C \longrightarrow C_{\sigma}\times C_{\tau}.
$$
Then, we have a homomorphism
\begin{equation}\label{Nf}
 N_f = (N_{f_1}, N_{f_2}): J(C) \longrightarrow J(C_{\sigma})\times J(C_{\tau}).
\end{equation}

The automorphisms $\sigma$ and $\tau$ induce the automorphisms of $J(C)$,
and we have natural isomorphisms:
$$
\begin{array}{l}
  {\rm H}^0(J(C), \Omega_{J(C)}^{1})\cong {\rm H}^0(C, \Omega_C^{1}) 
  \cong 
  {\rm H}^0(C, \Omega_C^{1})^{\langle \sigma^*\rangle} \oplus 
  {\rm H}^0(C, \Omega_C^{1})^{\langle \tau^*\rangle} \\
  \cong 
  {\rm H}^0(C_{\sigma}, \Omega_{C_{\sigma}}^{1}) \oplus {\rm H}^0(C_{\tau}, \Omega_{C_{\tau}}^{1}) \cong 
  {\rm H}^0(J(C_{\sigma}), \Omega_{J(C_{\sigma})}^{1}) \oplus 
  {\rm H}^0(J(C_{\tau}), \Omega_{J(C_{\tau})}^{1}).
\end{array}
$$
Therefore, $N_f$ is an isogeny. 
Note that 
$$
N_{f_1}\circ f_1^* =[2]_{J(C_{\sigma})}, \quad
N_{f_2}\circ f_2^* =[2]_{J(C_{\tau})}.
$$
By our construction, we have
$$
N_{f_1}\circ f_2^* =0, \quad
N_{f_2}\circ f_1^* =0.
$$
Therefore, we have
\begin{equation*}
N_{f}\circ f^* =[2]_{J(C_{\sigma})\times J(C_{\tau})}.
\end{equation*}
Dualizing the situation (\ref{Nf}), we have
$$
   f^* : J(C_{\sigma})\times J(C_{\tau}) \longrightarrow J(C).
$$

\begin{theorem}
\label{Thm:hyperelliptic}
Let $C$ be a hyperelliptic curve  
with an automorphism $\sigma$ of order 2, which is not the inversion.
We set $\tau = \sigma\circ \iota$ as above.
Then, the isogeny $N_f : J(C) \longrightarrow J(C_{\sigma})\times J(C_{\tau})$
is a decomposed Richelot isogeny.
\end{theorem}
\proof{
Since $\sigma$ induces 
an isomorphism from $J(C)$ to $J(C)$ and we may assume
that this isomorphism is an automorphism of $J(C)$,
we have a commutative diagram
$$
\begin{array}{ccc}
  J(C_{\sigma})\times J(C_{\tau})& \stackrel{id_{J(C_{\sigma})}\times \iota_{J(C_{\tau})}}{\longrightarrow}  & J(C_{\sigma})\times J(C_{\tau})\\
   f^*  \downarrow    &        & \quad \downarrow f^*  \\
    \quad J(C) & \stackrel{\sigma}{\longrightarrow }   & J(C) \\
  N_f   \downarrow &     & \quad \downarrow N_f\\
   J(C_{\sigma})\times J(C_{\tau})& \stackrel{id_{J(C_{\sigma})}\times \iota_{J(C_{\tau})}}{\longrightarrow} & J(C_{\sigma})\times J(C_{\tau}).
\end{array}  
$$
Since $\sigma^*(\Theta) = \Theta$, using Lemma \ref{decompose}, we have
$$
f^*(\Theta) \approx f_1^{*}(\Theta) \times J(C_{\tau})+ J(C_{\sigma}) \times f_2^*(\Theta).
$$
Therefore, by Lemma \ref{2Theta}, we see
$$
 f^*(\Theta) \approx 2(C_{\sigma} \times J(C_{\tau}))+ 2(J(C_{\sigma}) \times C_{\tau}).  
$$
Dualizing this situation, we have
$$
      N_f^*((C_{\sigma} \times J(C_{\tau}))+ (J(C_{\sigma}) \times C_{\tau}))\approx 2\Theta.
$$
This means that $N_f$ is a decomposed Richelot isogeny outgoing from $J(C)$.
\qed
}


\section{Generalized Howe curves}
\label{Sec:Howe}

Let $C_1$, $C_2$ be the nonsingular projective models of 
two hyperelliptic curves defined respectively by
$$
\begin{array}{l}
C_1 : y_1^2 = (x - a_1)(x - a_2)\ldots (x - a_r)(x - a_{r + 1})\ldots (x - a_{2g_1 + 2}),\\
C_2 : y_2^2 = (x - a_1)(x - a_2)\ldots (x - a_r)(x - b_{r + 1})\ldots (x - b_{2g_2 + 2})
\end{array}
$$
Here, $a_i\neq a_j$, $b_i \neq b_j$ for $i \neq j$, and $a_i \neq b_j$ for any $i, j$.
We assume $0< g_1 \leq g_2$. The genera of these curves are given by
$$
g(C_1) = g_1,~ g(C_2) = g_2.
$$
Let $\psi_1 : C_1\longrightarrow {\bf P}^1$ and  $\psi_2 : C_2 \longrightarrow {\bf P}^1$ 
be the hyperelliptic structures. We have $r$ common branch points 
of $\psi_1$ and $\psi_2$ ($0 \leq r \leq 2g_1 + 2)$).
We consider the fiber product $C_1 \times_{{\bf P}^1} C_2$:
$$
\begin{array}{ccc}
   C_1 \times_{{\bf P}^1} C_2 & \stackrel{\pi_2}{\longrightarrow} & C_2 \\
     \pi_1 \downarrow &        & \downarrow  \psi_2\\
      C_1     & \stackrel{\psi_1}{\longrightarrow}   & {\bf P}^1.
\end{array}
$$

We assume that there exists no isomorphism $\varphi : C_1 \longrightarrow C_2$
such that $\psi_2\circ \varphi = \psi_1$. Then, the curve $C_1 \times_{{\bf P}^1} C_2$
is irreducible.
We denote by $C$ the nonsingular projective model of $C_1 \times_{{\bf P}^1} C_2$,
and  we denote by $h : C \longrightarrow C_1 \times_{{\bf P}^1} C_2$ 
the resolution of singularities.
We call $C$ a {\em generalized Howe curve}. If $g_1= g_2 = 1$,
$C$ is called a {\em Howe curve}, which was originally defined 
in genus 4 by Howe \cite{How16} (see also Kudo--Harashita--Senda \cite{KHS20}, Oort \cite{O91} and van der Geer--van der Vlugt \cite{GV95}).  
The naming comes from  Kudo-Harashita-Senda, loc.\,cit.
In the case of genus 3, a Howe curve is nothing but a Ciani curve.
We set $f_i = \pi_i\circ h$ for $i=1,2$. Then, the degree of $f_i$ is 2.
We have the following proposition.

\begin{proposition}\label{Howe}
The genus $g(C)$ of $C$ is equal to $2(g_1 + g_2) + 1 - r$.
\end{proposition}
\proof{
Let $P \in {\bf P}^1$ be a common branch point of $\psi_1$ and $\psi_2$.
We can choose a coordinate $x$ on ${\bf A}^1 \subset {\bf P}^1$
such that $P$ is locally defined by $x = 0$. Then, the equation
of $C_1$ (resp. $C_2$) around $P$ is given by
$$
   y_1^2 = u_1 x\quad (\mbox{resp.}~y_2^2 = u_2 x).
$$
Here, $u_1$ and $u_2$ are units at $P$.
We denote by $\tilde{P}$ the point 
of the fiber product $C_1 \times_{{\bf P}^1} C_2$
over $P$. Then, around $\tilde{P}$ 
the fiber product $C_1 \times_{{\bf P}^1} C_2$
is defined by
$$
   y_1^2 = u_1 x, ~y_2^2 = u_2 x.
$$
Therefore, by eliminating $x$, the equation around $\tilde{P}$
is given by the equation $u_2y_1^2 = u_1y_2^2$. This means
that $\tilde{P}$ is a singular point with two branches.
Therefore, on $C$, $\tilde{P}$ splits into two nonsingular
points and  $P$ is not a branch point of $f_1$.

By the meaning of fiber product, the ramification points of $\psi_1$
whose images by $\psi_1$ are not branch points of $\psi_2$
are not branch points of $f_1$, and the points on
$C_1$ which are not 
ramification points of $\psi_1$ and whose images by $\psi_1$ are 
branch points of $\psi_2$ are branch points 
of $f_1$. Therefore, 
on the curve $C$, $f_1$ has $2(2g_2 + 2 - r)$ ramification points
of index 2. Applying the Hurwitz formula 
to the morphism $f_1: C \longrightarrow C_1$, we have
$$
   2(g(C) -1) = 2\cdot 2(g(C_1) -1) + 2(2g_2 + 2 - r).
$$
Therefore, we have $g(C) = 2(g_1 + g_2) + 1 - r$.
\qed
}
\vspace*{0.3cm}

We denote by $\iota_{C_1}$ (resp.~$\iota_{C_2}$) the hyperelliptic 
involution of $C_1$ (resp.~$C_2$). Then, these involutions lift
to automorphisms of $C$ as follows:
$$
\begin{array}{l}
 \sigma = \iota_{C_1}: y_1 \mapsto -y_1, y_2 \mapsto y_2, x \mapsto x,\\
 \tau = \iota_{C_2}: y_1 \mapsto y_1, y_2 \mapsto -y_2, x \mapsto x.
\end{array}
$$
Both $\sigma$ and $\tau$ are of order 2 and we have $\sigma \circ \tau = \tau \circ \sigma$.
Clearly, we have $C/\langle \sigma\rangle \cong C_2$ and $C/\langle \tau \rangle \cong C_1$.
We set $y_3 = y_1y_2/(x - a_1)(x - a_2)\cdots (x - a_r)$. Then, we have
a curve $C_3 = C/\langle \sigma\circ \tau \rangle$, which is given by the equation
$$
   y_3^2 = (x - a_{r + 1})\cdots (x - a_{2g_1 + 2})(x - b_{r + 1})\cdots (x - b_{2g_2 + 2}).
$$
Since the degree of the polynomial of right hand side is $2(g_1 + g_2) + 4 - 2r$,
the genus of the curve $C_3$ is given by
\begin{equation}\label{g(C_3)}
g(C_3) = g_1 + g_2 + 1 - r,
\end{equation}
and we have
\begin{equation}\label{g(C)}
g(C) = g(C_1) + g(C_2) + g(C_3).
\end{equation}
We have natural projections $f_i : C \longrightarrow C_i$ ($i = 1, 2, 3$).
We denote by $(J(C_i), \Theta_i)$ the Jacobian varieties of $C_i$ ($i = 1, 2, 3$).

\begin{theorem}\label{rationalHowe}
Under the notation above, assume $g(C) \geq 4$.
Then, the generalized Howe curve $C$ is hyperelliptic
if and only if $r = g_1 + g_2 + 1$, i.e., the curve $C_3$ is rational.
\end{theorem}
\proof{First, we note $r \leq g_1 + g_2 + 1$. Because if $r > g_1 + g_2 + 1$,
by $r \leq 2g_1 +2$ and $g_1\leq g_2$ we have $g_1 = g_2$ and $r = 2g_1 + 2$ and 
all the branch points on ${\bf P}^1$ of $C_1$ and $C_2$ coincide. Therefore, 
$C_1$ is isomorphic to $C_2$, and there exists an automorphism
$\varphi : C_1 \longrightarrow C_2$ such that $\psi_1 = \psi_2 \circ \varphi$.
Therefore, the fiber product $C_1 \times_{{\bf P}^1} C_2$
is reducible by the universality of fiber product, and
we already excluded this case.

If $r = g_1 + g_2 + 1$, then we have $g(C_3) = 0$ and we have a morphism 
$C \longrightarrow C_3$ of degree 2. Therefore, $C$ is hyperelliptic.
If $r < g_1 + g_2 + 1$, then we have $g(C_3) > 0$. 
Since we have $g(C) \geq 4$, by (\ref{g(C)}) 
there exists $C_i$ such that 
$g(C_i)\geq 2$ and we have a morphism $f_i: C \longrightarrow C_i$.
Since $f_i$ is separable, we have an injective homomorphism
\begin{equation}\label{injection}
   {\rm H}^0(C_i, \Omega_{C_i}^1) \longrightarrow {\rm H}^0(C, \Omega_{C}^1).
\end{equation}
Suppose $C$ is hyperelliptic. 
Note that $C_i$ is a hyperelliptic curve.
Since the hyperelliptic structure of $C$ (resp. $C_i$) is given by
${\rm H}^0(C, \Omega_{C}^1)$ (resp. ${\rm H}^0(C_i, \Omega_{C_i}^1)$)
we have the following commutative diagram by (\ref{injection}):
\begin{equation}\label{diagram}
\begin{array}{ccc}
    C  & \stackrel{f_i}{\longrightarrow} & C_i \\
    \downarrow &       &  \downarrow \\
   ~ {\bf P}^1 &\longrightarrow & ~ {\bf P}^1.
\end{array}
\end{equation}
We have a field extension $k(C)/k(x)$. This is a Galois extension 
and the Galois group is isomorphic to ${\bf Z}/2{\bf Z}\times {\bf Z}/2{\bf Z}$.
Therefore, we have 3 intermediate fields of the field extension $k(C)/k(x)$,
and the 3 intermediate fields are given by $k(C_i)$ ($i = 1, 2, 3$)
whose genera are greater than or equal to 1.
However, by the diagram (\ref{diagram}) we have one more intermediate field
$k({\bf P}^1)$, a contradiction.
\qed
}
\vspace*{0.3cm}

From Theorem \ref{rationalHowe}, the Howe curves of genus 4 which are constructed with 
$g_1=1$, $g_2 = 1$ and $r = 1$ are non-hyperelliptic, which is known
to Kudo-Harashita-Howe \cite[Lemma 2.1]{KHH20}.

\begin{theorem}\label{Thm:Howe}
Let $C$ be a generalized Howe curve defined as above. 
Then, $C$ has a decomposed Richelot isogeny
given by a natural isogeny 
$$
N_f : J(C) \longrightarrow J(C_1) \times J(C_2) \times J(C_3).
$$
\end{theorem}
\proof{
We have a homomorphism
$$
    N_f= (N_{f_1}, N_{f_2}, N_{f_3}) : J(C) \longrightarrow J(C_1) \times J(C_2) \times J(C_3).
$$
Since we have natural isomorphisms:
\begin{eqnarray*}
&& {\rm H}^0(J(C), \Omega_{J(C)}^{1}) \cong {\rm H}^0(C, \Omega_C^{1}) \\
&& \quad \quad \cong 
  {\rm H}^0(C, \Omega_C^{1})^{\langle \sigma^*\rangle} \oplus 
  {\rm H}^0(C, \Omega_C^{1})^{\langle \tau^*\rangle} \oplus
  {\rm H}^0(C, \Omega_C^{1})^{\langle \tau^*\circ \sigma^*\rangle}\\
&&  \quad \quad \cong 
  {\rm H}^0(C_{\sigma}, \Omega_{C_{\sigma}}^{1}) \oplus {\rm H}^0(C_{\tau}, \Omega_{C_{\tau}}^{1}) 
\oplus {\rm H}^0(C_{\sigma\circ \tau}, \Omega_{C_{\sigma \circ \tau}}^{1})\\
&&  \quad \quad  \cong {\rm H}^0(J(C_{1}), \Omega_{J(C_{1})}^{1}) \oplus {\rm H}^0(J(C_{2}), \Omega_{J(C_{2})}^{1}) 
\oplus {\rm H}^0(J(C_{3}), \Omega_{J(C_{3})}^{1}),
\end{eqnarray*}
we see that $N_f$ is an isogeny.
Then by a similar method to the one in Theorem \ref{Thm:hyperelliptic},
we have 
$$
2\Theta \approx {N_f}^*(\Theta_1 \times J(C_2)\times J(C_3) + 
J(C_1) \times \Theta_2 \times J(C_3) + J(C_1) \times J(C_2)\times \Theta_3)
$$
and $N_f$ is a decomposed Richelot isogeny.
\qed
}

\begin{remark}
\label{Rmk:hyperelliptic}
{\rm
Under the notation in Section \ref{Sec:Hyperelliptic}, we set $C_1 = C_{\sigma}$
and $C_2 = C_{\tau}$. Let
$\psi_{i} : C_{i}\longrightarrow {\bf P}^1$ ($i = 1, 2$)
be the hyperelliptic structures. Then, by the universality of fiber product
we have $C \cong C_{1}\times_{{\bf P}^1} C_2$. In this case, we have $\sigma\circ \tau = \iota$, 
and $C_3 = C/\langle \iota\rangle\cong {\bf P}^1$.
This means that the hyperelliptic curve $C$ is a generalized Howe curve 
which satisfies the condition
$r = g_1 + g_2 + 1$ in Theorem \ref{rationalHowe}, and that Theorem \ref{Thm:hyperelliptic} is a special case of Theorem \ref{Thm:Howe}.
}
\end{remark}

If the genus of the Howe curve $C$ is 3, then the converse of Theorem \ref{Thm:Howe}
holds
(cf.\,\cite[Theorem 6.3]{K}), that is, if there exists a completely
decomposed Richelot isogeny outgoing from the Jacobian variety $J(C)$,
then the curve $C$ is a Howe curve. However, if the genus of the curve $C$
is large, it seems to be difficult to formulate the converse.
Assume the target of the Richelot isogeny is decomposed into 3 principally
polarized abelian varieties as in Theorem \ref{Thm:Howe}. 
First, if the dimension of a principally polarized abelian variety is large, 
then it is not
necessarily a Jacobian variety. If the components of
the decomposition are all Jacobian varieties, the automorphism
of order 2 of the target does not necessarily determine
a good automorphism of $C$ if the curve is not hyperelliptic.
Note that for an automorphism $\sigma$ of order 2 of the non-hyperelliptic curve
$C$ of genus 3 the quotient curve $C/\langle \sigma\rangle$ is always
an elliptic curve (cf.\,Katsura \cite[Corollary 5.2]{K}).
Such facts work well in the case of genus 3. But we cannot expect
similar results in higher genus.


\section{Examples}
\label{Sec:Examples}
In this section,  we assume  the characteristic $p \neq 2$ and
give some concrete examples.


\begin{example}
\label{Ex:NonHyperelliptic}
\rm{
We consider the nonsingular complete model $C$ of a curve defined by the equation
$$
   x^4 + y^4 + x^2y^2 + 1 = 0.
$$
The genus of this curve is 3 and $C$ has automorphisms $\sigma$, $\tau$ of order 2 given by
$$
\begin{array}{l}
   \sigma: x \mapsto -x,~ y \mapsto y\\
   \tau : x \mapsto x,~ y \mapsto -y.
\end{array}
$$
We set $u = \sqrt[4]{3} y/\sqrt{2}$ and $v = x^2 + (y^2/2)$ (resp. 
$u = \sqrt[4]{3} x/\sqrt{2}$ and $v = y^2 + (x^2/2)$). 
Then, $u, v$ are invariant under the action of the group $\langle \sigma \rangle$
(resp.~the group $\langle \tau \rangle$) and the quotient curve 
$E_{\sigma} = C/\langle \sigma \rangle$ (resp.~$E_{\tau} = C/\langle \tau \rangle$)
is an elliptic curve defined by the equation
$$
     v^2 +  u^4 + 1 = 0.
$$
Since the group $G =\langle \sigma, \tau \rangle \cong {\bf Z}/2{\bf Z} \times {\bf Z}/2{\bf Z}$
acts on $C$ and we have $C/G \cong {\bf P}^1$, we see that the original curve $C$ is
a non-hyperelliptic Howe curve given by the fiber product $E_{\sigma}\times_{{\bf P}^1}E_{\tau}$.
Since $u = y/x$ and $v = 1/x^2$ are invariant under the action of
the group $\langle \sigma \circ \tau \rangle$, we have the third elliptic curve 
$E_{\sigma\circ \tau} = C/\langle \sigma \circ \tau\rangle$ defined by
$$
     v^2 + u^4 + u^2 + 1 = 0.
$$
We have a natural morphism $C \longrightarrow E_{\sigma}\times E_{\tau}\times E_{\sigma\circ \tau}$,
and this morphism induces a completely decomposed Richelot isogeny
$$
    J(C) \longrightarrow E_{\sigma}\times E_{\tau}\times E_{\sigma\circ \tau}.
$$
}
\end{example}


\begin{example}
\label{Ex:Hyperelliptic}
\rm{
We consider the nonsinglar complete model $C$ of a hyperelliptic curve defined by the equation
$$
   y^2 = x^8 + x^4 + 1.
$$
The genus of this curve is 3 and $C$ has automorphisms $\sigma$, $\iota$ of order 2 given by
$$
\begin{array}{l}
   \sigma: x \mapsto -x,~ y \mapsto y,\\
   \iota : x \mapsto x,~ y \mapsto -y.
\end{array}
$$
The automorphism $\iota$ is a hyperelliptic involution.
We set $u = x^2$ and $v = y$ (resp.~$u = x^2$ and $v = xy$). Then, $u, v$ are invariant under the action of 
the group $\langle \sigma \rangle$
(resp.~the group $\langle \sigma \circ \iota \rangle$) and the quotient curve 
$E_{\sigma} = C/\langle \sigma \rangle$ 
(resp.~$C_{\sigma \circ \iota} = C/\langle \sigma \circ \iota \rangle$)
is a curve defined by the equation
$$
     v^2 = u^4 + u^2 + 1 \quad (\mbox{resp.}~ v^2 = u(u^4 + u^2 + 1)).
$$
$E_{\sigma}$ is an elliptic curve and $C_{\sigma \circ \iota}$ is a curve of genus 2.
We have a natural morphism $C \longrightarrow E_{\sigma}\times C_{\sigma \circ \iota}$,
and this morphism induces a decomposed Richelot isogeny
$$
          J(C) \longrightarrow E_{\sigma}\times J(C_{\sigma \circ \iota}).
$$

On the other hand, we consider the following automorphism:
$$ 
    \tau : x \mapsto 1/x,~ y \mapsto y/x^4.
$$
We set $u = x + (1/x)$ and $v = y/x^2$ (resp.~$u = x - (1/x)$ and $v = y/x^2$). 
Then, $u$ and $v$ are invariant under the action
of the group $\langle \tau \rangle$ (resp.~the group $\langle \sigma \circ \tau \rangle$) and the quotient curve $E_{\tau} = C/\langle \tau \rangle$ 
(resp. $E_{\sigma\circ \tau} = C/\langle \sigma\circ \tau \rangle$) is an elliptic curve and
given by the equation
$$
     v^2 = u^4 - 4u^2 + 3\quad (\mbox{resp.}~v^2 = u^4 + 4u^2 + 3).
$$
Since the group 
$G =\langle \tau, \sigma \circ \tau \rangle \cong {\bf Z}/2{\bf Z} \times {\bf Z}/2{\bf Z}$
acts on $C$ and we have $C/G \cong {\bf P}^1$, we see that the original curve $C$ is
a hyperelliptic Howe curve given 
by the fiber product $E_{\tau}\times_{{\bf P}^1}E_{\sigma \circ \tau}$.
We have 
a natural morphism $C \longrightarrow E_{\sigma} \times E_{\tau}\times E_{\sigma \circ \tau}$
and this induces a completely decomposed Richelot isogeny
$$
J(C) \longrightarrow E_{\sigma}\times E_{\tau}\times E_{\sigma\circ \tau}.
$$
}
\end{example}


\begin{example}\label{Fermat}
\rm{
We consider the nonsingular complete model $\FermatC$ of a curve defined by the equation
\begin{eqnarray*}
\label{Eq:Fermat}
   \FermatC: \ x^4 + y^4 + 1 = 0.
\end{eqnarray*}
The genus of this curve is 3 and $\FermatC$ has automorphisms $\sigma$, $\tau$ of order 2 given by
$$
\begin{array}{l}
   \sigma: x \mapsto -x,~ y \mapsto y,\\
   \tau : x \mapsto x,~ y \mapsto -y.
\end{array}
$$
We set $u = y$ and $v =x^2$ (resp.~$u = x$ and $v = y^2$). 
Then, $u, v$ are invariant under the action of the group $\langle \sigma \rangle$
(resp.~the group $\langle \tau \rangle$) and the quotient curve 
$E_{\sigma} = \FermatC/\langle \sigma \rangle$ (resp.~$E_{\tau} = \FermatC/\langle \tau \rangle$)
is an elliptic curve defined by the equation
\begin{equation}\label{elliptic}
     v^2 +  u^4 + 1 = 0.
\end{equation}
Since the group $G =\langle \sigma, \tau \rangle \cong {\bf Z}/2{\bf Z} \times {\bf Z}/2{\bf Z}$
acts on $\FermatC$ and we have $\FermatC/G \cong {\bf P}^1$, we see that the original curve $\FermatC$ is
a non-hyperelliptic Howe curve given by the fiber product $E_{\sigma}\times_{{\bf P}^1}E_{\tau}$.
Since $u = y/x$ and $v = 1/x^2$ are invariant under the action of
the group $\langle \sigma \circ \tau \rangle$, we have the third elliptic curve 
$E_{\sigma\circ \tau} = \FermatC/\langle \sigma \circ \tau\rangle$ defined 
by the equation (\ref{elliptic}).
We have a natural morphism $\FermatC \longrightarrow E_{\sigma}\times E_{\tau}\times E_{\sigma\circ \tau}$,
and this morphism induces a completely decomposed Richelot isogeny
$$
    J(\FermatC) \longrightarrow E_{\sigma}\times E_{\tau}\times E_{\sigma\circ \tau}.
$$
Since the elliptic curve defined by the equation (\ref{elliptic}) has automorphism
of order 4, it is isomorphic to the elliptic curve $E_0$ defined by
\begin{eqnarray*}
\label{Eq:E0}
  E_0: \  y^2 = x^3  - x
\end{eqnarray*}
over an algebraically closed field $k$. 
$E_0$ is supersingular if and only if $p \equiv 3~(\mbox{mod}~ 4)$.
Since our Richelot isogeny is separable, we see that the curve
$\FermatC$ is a superspecial non-hyperelliptic Howe curve if $p \equiv 3~(\mbox{mod}~ 4)$, that is, 
the Jacobian variety $J(\FermatC)$
is isomorphic to a product of three supersingular elliptic curves.
The elliptic curve $E_0$ has an automorphism $\rho$ of order 4 defined by
$$
 \rho : x \mapsto -x ,~y \mapsto iy.
$$
Here, $i$ is a primitive fourth root of unity. We denote by 
$F$ the Frobenius morphism of $E_0$. We note that if $p \equiv 3~(\mbox{mod}~ 4)$,
the endomorphism ring of $E_0$ is given by
\begin{eqnarray*}
\label{Eq:EndE0}
{\rm End}(E_0) = {\bf Z}+ {\bf Z}\rho + {\bf Z}(1 + F)/2 + {\bf Z}\rho (1 + F)/2
\end{eqnarray*}
(cf.\,Katsura \cite{KT87}).
}
\end{example}

\begin{remark}
{\rm
In the case of curves of genus 3, decomplosed Richelot isogenies  are 
studied in
Howe-Lepr\'evost-Poonen \cite{HLP} and Katsura \cite{K} in detail.
}
\end{remark}

We give two more examples of higher genera. By this method, we can construct
many superspecial curves (see also Kudo-Harashita-Howe \cite{KHH20}).


\begin{example}
\label{Ex:NonHyperellipticLargeGenus}
{\rm
Assume the characteristic $p > 2$.
We consider two elliptic curves
$$
    C_1: \y_1^2 = x^4 + 1, \quad C_2 : \y_2^2= x^4 - 1.
$$
They are isomorphic to each other and 
supersingular if and only if $p \equiv 3~({\rm mod}~4)$.
Let $C$ be the generalized Howe curve which is birational to $C_1\times_{{\bf P}^1}C_2$.
Then we have a Galois extension $k(C)/k(x)$ with the Galois group 
$\cong {\bf Z}/2{\bf Z} \times {\bf Z}/2{\bf Z}$ and we have $k(C) \cong k(x, \y_1, \y_2)$.
In this case we have $r = 0$, and by the formula (\ref{g(C_3)}), we have $g(C_3) = 3$.
We set $\y_3 = \y_1\y_2$. Then we have
$$
       \y_3^2 = x^8 - 1
$$
which is the equation for the curve $C_3$ of genus 3 and we have three intermediate
field $k(C_1)$, $k(C_2)$ and $k(C_3)$ of the field extension $k(C)/k(x)$.
By the calculation of the Cartier operator, we can easily show that
$C_3$ is superspecial if and only if $p \equiv 7~({\rm mod}~8)$.
Now, we set $y = \y_1 + \y_2$. Then, we have $\y_1\y_2 = y^2/2 -x^4$.
Using this equation, we know that the curve $C$
is the nonsingular model of the curve defined by the following equation:
$$
     y^4 = 4x^4y^2 - 4.
$$
Since we have $r = 0$, 
$C$ is non-hyperelliptic by Theorem \ref{rationalHowe}.
Since we have a Richelot isogeny 
$J(C) \longrightarrow J(C_1) \times J(C_2) \times J(C_3)$ and 
the Richelot isogeny is separable, the curve $C$ is superspecial
if and only if $p \equiv 7~({\rm mod}~8)$.

Incidentally, the three automorphims of order 2 of the curve $C$ are given by
$$
\sigma : x \mapsto x, ~y \mapsto 2/y,\quad  \tau: x \mapsto x,~y \mapsto - 2/y,
$$
and $\sigma \circ \tau$.
}
\end{example}


\begin{example}
\label{Ex:HyperellipticLargeGenus}
{\rm 
Assume the characteristic $p \geq 7$.
We consider two curves of genus 2:
$$
    C_1: \y_1^2 = x^5 + 1, \quad C_2 : \y_2^2= x^6 + x.
$$
By the isomorphism
$$
  x \mapsto 1/x,~ \y_2 \mapsto \y_1/x^3
$$
they are isomorphic to each other.
Moreover, they are supersingular if and only if $p \equiv 4~({\rm mod}~5)$
(cf.\,Ibukiyama-Katsura-Oort \cite{IKO86}).
Let $C$ be the generalized Howe curve which is birational to $C_1\times_{{\bf P}^1}C_2$.
Then we have a Galois extension $k(C)/k(x)$ with the Galois group 
$\cong {\bf Z}/2{\bf Z} \times {\bf Z}/2{\bf Z}$ and we have $k(C) \cong k(x, \y_1, \y_2)$.
In this case we have $r = 5$, and by the formula (\ref{g(C_3)}), we have $g(C_3) = 0$.
Therefore, $C$ is hyperelliptic and we have $g(C) = 4$ by (\ref{g(C)}).
We set $\y_3 = \y_2/\y_1$. Then we have the equation of $C_3$:
\begin{equation}\label{y}
       \y_3^2 = x
\end{equation}
and we have three intermediate
fields $k(C_1)$, $k(C_2)$ and $k(C_3)$ of the field extension $k(C)/k(x)$.
Now, we set $y = \y_1 + \y_2$.  Then we have
\begin{equation}\label{yv}
     y^2 = (x^5 + 1)(1 + x + 2\y_3).
\end{equation}
Using the equation   (\ref{y}), and putting $z = y/(1 + \y_3)$, we have
the equation of the curve $C$ over $C_3$:
$$
       z^2 = \y_3^{10} + 1.
$$
Using the equations (\ref{yv}) and (\ref{y}), we have the equation of the curve $C$:
$$
  y^4 -2(x^5 + 1)(x + 1)y^2 + (x^5+ 1)^2(x-1)^2= 0.
$$
Since we have a Richelot isogeny 
$J(C) \longrightarrow J(C_1) \times J(C_2)$ and 
the Richelot isogeny is separable, the curve $C$ is superspecial
if and only if $p \equiv 4~({\rm mod}~5)$.

Incidentally, the three automorphims of order 2 of the curve $C$ are given by
$$
\begin{array}{l}
\sigma : x \mapsto x, ~y \mapsto -(x-1)(x^5 + 1)/y,\\
\tau: x \mapsto x,~y \mapsto (x-1)(x^5 + 1)/y,
\end{array}
$$
and $\sigma \circ \tau$.
}
\end{example}

\bibliographystyle{spmpsci}
\bibliography{StPauli_2023}

\begin{thebibliography}{10}
\providecommand{\url}[1]{{#1}}
\providecommand{\urlprefix}{URL }
\expandafter\ifx\csname urlstyle\endcsname\relax
  \providecommand{\doi}[1]{DOI~\discretionary{}{}{}#1}\else
  \providecommand{\doi}{DOI~\discretionary{}{}{}\begingroup
  \urlstyle{rm}\Url}\fi

\bibitem{BL}
Birkenhake, C., Lange, H.: Complex Abelian Varieties, 2nd edn.
\newblock Springer-Verlag Berlin Heidelberg (2004)

\bibitem{CD23}
Castryck, W., Decru, T.: An efficient key recovery attack on {SIDH}.
\newblock In: {EUROCRYPT} 2023, Part {V}, \emph{LNCS}, vol. 14008, pp.
  423--447. Springer (2023)

\bibitem{CDS}
Castryck, W., Decru, T., Smith, B.: Hash functions from superspecial genus-2
  curves using {R}ichelot isogenies.
\newblock J.~Math.~Crypt. \textbf{14}(1), 268--292 (2020)

\bibitem{ChaGorLau09}
Charles, D., Goren, E., Lauter, K.: Families of {Ramanujan} graphs and
  quaternion algebras.
\newblock In: Groups and Symmetries: From Neolithic Scots to John McKay, pp.
  53--80 (2009)

\bibitem{CS20}
Costello, C., Smith, B.: The supersingular isogeny problem in genus 2 and
  beyond.
\newblock In: PQCrypto 2020, \emph{LNCS}, vol. 12100, pp. 151--168. Springer
  (2020)

\bibitem{FJP14}
{De Feo}, L., Jao, D., Pl{\^{u}}t, J.: Towards quantum-resistant cryptosystems
  from supersingular elliptic curve isogenies.
\newblock J.~Math.~Crypt. \textbf{8}(3), 209--247 (2014)

\bibitem{DK23}
Decru, T., Kunzweiler, S.: Efficient computation of $(3^n,3^n)$-isogenies.
\newblock {IACR} Cryptol. ePrint Arch. 2023/376  (2023).
\newblock \urlprefix\url{https://eprint.iacr.org/2023/376}.
\newblock To appear at AFRICACRYPT 2023

\bibitem{FS21a}
Florit, E., Smith, B.: {Automorphisms and isogeny graphs of abelian varieties,
  with applications to the superspecial Richelot isogeny graph}.
\newblock In: {Arithmetic, Geometry, Cryptography, and Coding Theory 2021},
  \emph{Contemporary Mathematics}, vol. 779. {American Mathematical Society}
  (2021)

\bibitem{FS21b}
Florit, E., Smith, B.: {An atlas of the Richelot isogeny graph}.
\newblock In: Theory and Applications of Supersingular Curves and Supersingular
  Abelian Varieties, \emph{{RIMS K{\^o}ky{\^u}roku Bessatsu}}, vol. B90, pp.
  195--219. {RIMS, Kyoto University} (2022)

\bibitem{FT19}
Flynn, E.V., Ti, Y.B.: Genus two isogeny cryptography.
\newblock In: PQCrypto 2019, \emph{LNCS}, vol. 11505, pp. 286--306. Springer
  (2019)

\bibitem{GV95}
van~der Geer, G., van~der Vlugt, M.: On the existence of supersingular curves
  of given genus.
\newblock J. Reine Angew. Math. \textbf{458}, 53--61 (1995)

\bibitem{How16}
Howe, E.W.: Quickly constructing curves of genus 4 with many points.
\newblock In: D.~Kohel, I.~Shparlinski (eds.) {Frobenius
  distributions:Lang-Trotter and Sato-Tate conjectures}, \emph{Comtemp. Math.},
  vol. 663, pp. 149--173. Amer. Math. Soc., Providence, RI (2016)

\bibitem{HLP}
Howe, E.W., Lepr\'evost, F., Poonen, B.: Large torsion subgroups of split
  jacobians of curves of genus two or three.
\newblock Forum Math. \textbf{12}, 315--364 (2000)

\bibitem{IKO86}
Ibukiyama, T., Katsura, T., Oort, F.: {Supersingular curves of genus two and
  class numbers}.
\newblock Compositio Math. \textbf{57}, 127--152 (1986)

\bibitem{JZ}
Jordan, B.W., Zaytman, Y.: Isogeny graphs of superspecial abelian varieties and
  generalized {B}randt matrices.
\newblock ArXiv abs/2005.09031  (2020).
\newblock \urlprefix\url{https://arxiv.org/abs/2005.09031}

\bibitem{KT87}
Katsura, T.: Generarized {K}ummer surfaces and their unirationality in
  characteristic $p$.
\newblock J. Fac. Sci. Univ. Tokyo Sect. IA Math. \textbf{34}, 1--41 (1987)

\bibitem{K}
Katsura, T.: Decomposed {R}ichelot isogenies of {Jacobian} varieties of curves
  of genus 3.
\newblock J. Algebra \textbf{588}, 129--147 (2021)

\bibitem{KT20}
Katsura, T., Takashima, K.: Counting {R}ichelot isogenies between superspecial
  abelian surfaces.
\newblock In: ANTS 2020, \emph{The Open Book Series}, vol.~4, pp. 283--300.
  Mathematical Sciences Publishers (2020)

\bibitem{KH22}
Kudo, M., Harashita, S.: Algorithmic study of superspecial hyperelliptic curves
  over finite fields.
\newblock Commentarii Mathematici Universitatis Sancti Pauli \textbf{70},
  49--64 (2022)

\bibitem{KHH20}
Kudo, M., Harashita, S., Howe, E.: Algorithms to enumerate superspecial {Howe}
  curves of genus four.
\newblock In: ANTS 2020, \emph{The Open Book Series}, vol.~4, pp. 301--316.
  Mathematical Sciences Publishers (2020)

\bibitem{KHS20}
Kudo, M., Harashita, S., Senda, H.: The existence of supersingular curves of
  genus 4 in arbitrary characteristic.
\newblock Research in Number Theory \textbf{6}(4), 44 (2020)

\bibitem{Kum15}
Kumar, A.: Hilbert modular surfaces for square discriminants and elliptic
  subfields of genus 2 function fields.
\newblock Research in the Mathematical Sciences \textbf{2}(1), 1--46 (2015)

\bibitem{MMPPW23}
Maino, L., Martindale, C., Panny, L., Pope, G., Wesolowski, B.: A direct key
  recovery attack on {SIDH}.
\newblock In: {EUROCRYPT} 2023, Part {V}, \emph{LNCS}, vol. 14008, pp.
  448--471. Springer (2023)

\bibitem{MK22b}
Moriya, T., Kudo, M.: Computation of superspecial howe curves of genus 3.
  {Magma codes}.
\newblock
  \urlprefix\url{https://sites.google.com/view/m-kudo-official-website/english/code/genus3v4}.
\newblock Last accessed at 2023/03/20

\bibitem{MK22}
Moriya, T., Kudo, M.: Some explicit arithmetics on curves of genus three and
  their applications.
\newblock CoRR abs/2209.02926  (2022).
\newblock \urlprefix\url{https://doi.org/10.48550/arXiv.2209.02926}

\bibitem{M}
Mumford, D.: Abelian Varieties.
\newblock Oxford Univ. Press (1970)

\bibitem{O91}
Oort, F.: Hyperelliptic supersingular curves.
\newblock In: Arithmetic algebraic geometry (Texel, 1989), G. van der Geer, F.
  Oort, and J. Steenbrink, eds., \emph{Progr. Math.}, vol.~89, pp. 247--284.
  Birkh\"auser Boston, Boston, MA (1991)

\bibitem{Rob23}
Robert, D.: Breaking {SIDH} in polynomial time.
\newblock In: {EUROCRYPT} 2023, Part {V}, \emph{LNCS}, vol. 14008, pp.
  472--503. Springer (2023)

\bibitem{SCF22}
Santos, M.C., Costello, C., Frengley, S.: An algorithm for efficient detection
  of {$(N, N)$}-splittings and its application to the isogeny problem in
  dimension 2.
\newblock {IACR} Cryptol. ePrint Arch. 2022/1736  (2022).
\newblock \urlprefix\url{https://eprint.iacr.org/2022/1736}

\bibitem{Tak17}
Takashima, K.: Efficient algorithms for isogeny sequences and their
  cryptographic applications.
\newblock In: Mathematical Modelling for Next-Generation Cryptography: {CREST}
  Crypto-Math Project, pp. 97--114. {Springer} (2017)

\end{thebibliography}

\end{document}